\documentclass{article}

\usepackage{amssymb}
\usepackage{amsmath}

\usepackage[dvips]{graphicx}

\begin{document}


\begin{center}
\textbf{POSITIVE LINEAR FUNCTIONALS ON BP*-ALGEBRAS }
\end{center}
 
\noindent \textbf{}

\begin{center}
M. EL AZHARI 
\end{center}

\noindent \textbf{}

\noindent \textbf{Abstract.} Let $A$ be a BP*-algebra with identity $e,P_{1}(A)$ be the set of all positive linear functionals f on $A$ such that $f(e)=1,$ and let $M_{s}(A)$ be the set of all nonzero hermitian multiplicative linear functionals on $A.$ We prove that $M_{s}(A)$ is the set of extreme points of $P_{1}(A).$ We also prove that, if $M_{s}(A)$ is equicontinuous, then every positive linear functional      
on $A$ is continuous. Finally, we give an example of a BP*-algebra whose topological dual is not included in the vector space generated by  $P_{1}(A),$ which gives a negative answer to a question posed by M. A. Hennings.

\noindent \textbf{}


\noindent \textbf{} 


\noindent \textbf{}
 
\noindent \textbf{I. Preliminaries }                                                                        

\noindent \textbf{}
 
Let $A$ be a topological vector space, we denote by $A^{\ast}$ the algebraic dual of $A$, and by $ A^{\prime}$ the topoloçgical dual of $A.$ Let $A$ be an algebra and $ x\in A,$ we denote by $ sp_{A}(x)$
the spectrum of $x,$ by $ \rho_{A}(x)=\sup \lbrace\vert\lambda\vert:\lambda\in sp_{A}(x)\rbrace$ the spectral radius of $x,$ and by $R(A)$ the radical of $A.$ A topological algebra is an algebra (over the complex field) which is also a Hausdorff topological vector space such that the multiplication is separately continuous. A locally convex algebra is a topological algebra whose topology is locally convex. Let $A$ be a locally convex *-algebra. A linear functional $f$ on $A$ is called positive if $ f(x^{\ast}x)\geq 0 $ for all $ x\in A.$ A linear functional $f$ on $A$ is said to be hermitian if $ f(x^{\ast})=\overline{f(x)}$ for all  $ x\in A.$ Let $A$ be a commutative locally convex *-algebra with continuous involution and identity $e.$ Let $\beta$ be the family of all closed, bounded, absolutely convex, containing $e,$ hermitian, idempotent subsets of $A.$ For any $B\in\beta,$ the linear span $A_{B}$   
of $B,$ is a *-subalgebra of $A$ which, when equipped with the Minkowski norm $\Vert .\Vert_{B},$ becomes a normed *-algebra with identity and isometric involution. If $A=\cup\lbrace A_{B}, B\in\beta\rbrace$ and $(A_{B},\Vert .\Vert_{B})$ is complete for all $B\in\beta,$ we say that $A$ is a BP*-algebra. We denote by $P(A)$ the set of all positive linear functionals on $A,$ by $P_{1}(A)$ the subset of $P(A)$ consisting of those positive linear functionals $f$ such that $f(e)=1,$ by $M(A)$ the set of all nonzero multiplicative linear functionals on $A,$ and by $M_{s}(A)$ the set of all nonzero hermitian multiplicative linear functionals on $A.$ A BP*-algebra $A$ is called symmetric if $(e+x^{\ast}x)^{-1}$ exists for all $x\in A.$
 
\noindent \textbf{}

\noindent \textbf{II. Characterization of extreme points of $P_{1}(A)$}

\noindent \textbf{}

Let $A$ be a BP*-algebra, $A=\cup\lbrace A_{B}, B\in\beta\rbrace.$ In [5], T.Husain and S.A.Warsi tried to characterize the extreme points of $P_{1}(A),$ their method is essentially based on the fact that    $P_{1}(A)$ is the the projective limit of $P_{1}(A_{B}), B\in\beta,$ this method gave no result in the general case. Using the known method of Banach algebras, we characterize the extreme points of $P_{1}(A).$

\noindent \textbf{}

\noindent \textbf{Lemma 1.} ([6; Theorem 4.4.12]) Let $H$ be a Hilbert space and $L$ any *-subalgebra of  $B(H).$ Then $L$ is irreductible on $H$ if and only if the centralizer of $L$ in $B(H)$ contains only scalar multiples of the identity operator.

\noindent \textbf{}
 
\noindent \textbf {Proof.} See [6].

\noindent \textbf{}

\noindent \textbf{Theorem 2.} Let $A$ be a BP*-algebra, then the extreme points of $P_{1}(A)$ are exactly the nonzero hermitian multiplicative linear functionals on $A.$

\noindent \textbf{}
 
\noindent \textbf{Proof.} Let $\chi\in M_{s}(A)$ and $x\in A, \chi(x^{\ast}x)=\overline{\chi(x)}\chi(x)=\vert\chi(x)\vert^{2}\geq 0,$ therefore $M_{s}(A)\subset P_{1}(A).$ Let $f\in P_{1}(A), f$ is hermitian and admissible [4; Theorem 5.2], so we can associate to $f$ a cyclic representation $T_{f}$ on a Hilbert space $H_{f}$ [6; Theorem 4.5.4]. Let $f$ be an extreme point of $P_{1}(A), T_{f}$ is an irreductible representation [5; Theorem 5.1]. By Lemma 1, the centralizer $Z$ of $\lbrace T_{f}x, x\in A\rbrace$ in $B(H_{f})$ contains only scalar multiples of the identity operator. Since A is commutative, $\lbrace T_{f}x, x\in A\rbrace\subset Z,$ so for each $x\in A$ there exists $\chi(x)\in C$ such that $T_{f}x=\chi(x)I.$ It is easily shown that the map $x\rightarrow\chi(x)$ from     $A$ to $C$ is a nonzero hermitian multiplicative linear functional on $A.$ Let $h$ be a cyclic vector of  $H_{f}$ such that $f(x)=\langle T_{f}x(h),h\rangle$ for all $x\in A,$ then  $f(x)=\chi(x)\Vert h\Vert^{2},$ thus $\Vert h\Vert^{2}=1$ since $f(e)=\chi(e)=1,$ hence $f=\chi\in M_{s}(A).$ Conversely, let $\chi\in M_{s}(A)$ and let $(T_{\chi}, H_{\chi})$ be the cyclic representation associated with $\chi.$ We have dim $H_{\chi}=$ dim $A/Ker \chi = 1,$ therefore $T_{\chi}$ is irreductible, thus $\chi$ is an extreme point of     $P_{1}(A)$ [5; Theorem 5.1].
  
\noindent \textbf{}

\noindent \textbf{III. On the symmetric spectral radius of a BP*-algebra}

\noindent \textbf{}

Let $A$ be a BP*-algebra, we denote by $\rho_{s}$ the symmetric spectral radius of $A$ defined by $\rho_{s}(x)=\sup \lbrace\vert\chi(x)\vert: \chi\in M_{s}(A)\rbrace$ for all $x\in A.$ We have $\rho_{s}(x)\leq\rho_{A}(x)<\infty$ for all $x\in A.$

\noindent \textbf{}

\noindent \textbf{Theorem 3}. Let $A$ be a BP*-algebra, then $\rho_{s}(x)=\sup\lbrace f(x^{\ast}x): f\in P_{1}(A)\rbrace^{\frac{1}{2}}$ for all $x\in A.$

\noindent \textbf{}
 
\noindent \textbf{Proof.} Put $\vert x\vert=\sup\lbrace f(x^{\ast}x): f\in P_{1}(A)\rbrace^{\frac{1}{2}}$ for all $x\in A.$ Let $\chi\in M_{s}(A)\subset P_{1}(A)$ and $x\in A, \vert\chi(x)\vert^{2}=\chi(x^{\ast}x)\leq \vert x\vert^{2},$ thus $\vert\chi(x)\vert\leq\vert x\vert,$ hence $\rho_{s}(x)\leq\vert x\vert.$ Conversely, let $f\in P_{1}(A),$ by Theorem 2 and the Krein-Milman Theorem, there exists a net $(f_{\alpha})_{\alpha\in\Lambda}$ in conv $M_{s}(A)$ such that $f_{\alpha}\rightarrow f$
in the weak topology $\sigma(A^{\ast},A).$ For each $\alpha\in\Lambda,$ there exist $\lambda_{1},...,\lambda_{n}\in R^{+},\lambda_{1}+...+\lambda_{n}=1,$ and $\chi_{1},...,\chi_{n}\in M_{s}(A)$ such that     $f_{\alpha}=\lambda_{1}\chi_{1}+...+\lambda_{n}\chi_{n}.$ Let $x\in A, f_{\alpha}(x^{\ast}x)\leq \lambda_{1}\vert \chi_{1}(x)\vert^{2}+...+\lambda_{n}\vert \chi_{n}(x)\vert^{2}\leq \rho_{s}(x)^{2},$
since $ f_{\alpha}(x^{\ast}x)\rightarrow f(x^{\ast}x),$ it follows that $f(x^{\ast}x)\leq \rho_{s}(x)^{2},$ therefore $\sup\lbrace f(x^{\ast}x): f\in P_{1}(A)\rbrace\leq\rho_{s}(x)^{2},$ i.e.                $\vert x\vert\leq\rho_{s}(x).$  

\noindent \textbf{}

\noindent \textbf{Corollary 4.} Let $A$ be a symmetric BP*-algebra, then $R(A)=\lbrace x\in A: f(x^{\ast}x)=0$ for all $f\in P_{1}(A)\rbrace.$  
 
\noindent \textbf{}

\noindent \textbf{Proof.} Each nonzero multiplicative linear functional on $A$ is hermitian, therefore
$R(A)=\lbrace x\in A: \chi(x)=0$ for all $\chi\in M_{s}(A)\rbrace=\lbrace x\in A: \rho_{s}(x)=0\rbrace.$ 

\noindent \textbf{}

\noindent \textbf{Theorem 5.} Let $A$ be a BP*-algebra and $f\in P(A),$ then $\vert f(x)\vert\leq f(e)\rho_{s}(x)$ for all $x\in A.$

\noindent \textbf{}

\noindent \textbf{Proof.} Let $f\in P(A).$ If $f(e)=0,$ by the Schwarz inequality, $ \vert f(x)  \vert^{2}\leq f(e)f(x^{\ast}x)$ for all $x\in A,$ then $f$ vanishes on $A.$ If $f(e)\neq 0,$ put $f_{1}=f(e)^{-1}f,$ it is clear that $f_{1}\in P_{1}(A).$ Let $x\in A, f_{1}(x^{\ast}x)\leq\rho_{s}(x)^{2}$ by Theorem 3, thus $\vert f_{1}(x)\vert^{2}\leq f_{1}(x^{\ast}x)\leq \rho_{s}(x)^{2},$ therefore
$\vert f_{1}(x)\vert\leq\rho_{s}(x),$ i.e. $ \vert f(x)\vert\leq f(e)\rho_{s}(x).$

\noindent \textbf{}

\noindent \textbf{Corollary 6.} Let $A$ be a BP*-algebra. If $M_{s}(A)$ is equicontinuous, then every positive linear functional on $A$ is continuous.

\noindent \textbf{}

\noindent \textbf{Proof.} If $M_{s}(A)$ is equicontinuous, then $\rho_{s}$ is continuous.   

\noindent \textbf{} 

\noindent \textbf{Corollary 7.} Let $A$ be a barrelled BP*-algebra such that every nonzero hermitian multiplicative linear functional on $A$ is continuous. Then every positive linear functional on $A$ is continuous.
 
\noindent \textbf{}

\noindent \textbf{Proof} $M_{s}(A)$ is bounded in the weak topology $\sigma (A',A).$ Since $A$ is barrelled, $M_{s}(A)$ is equicontinuous by the Banach-Steinhaus Theorem.

\noindent \textbf{} 

\noindent \textbf{Remarks.} \\
1. Corollary 6 is an improvement of [4; Theorem 5.7].\\
2. Corollary 7 is a partial answer to the following problem [1]: is every positive linear functional on a barrelled BP*-algebra continuous ?
 
\noindent \textbf{} 

\noindent \textbf{IV. On a question of M. A. Hennings}  
 
\noindent \textbf{}  

Let $A$ be a BP*-algebra. In [3; question E], M. A. Hennings posed the following question:
do the positive linear functionals on $A$ spans the topological dual of $A$ ? \\
\noindent \textbf{}M. A. Hennings [3, p.290 (lines 28-29)] stated that the general problem of whether this question is true or not is still open. Here we give a negative answer to this question, we first have the following proposition:

\noindent \textbf{}

\noindent \textbf{Proposition 8.} Let $A$ be a symmetric BP*-algebra. If $A'\subset span_{A^{\ast}}P(A),$ then $A$ is semisimple.

\noindent \textbf{}

\noindent \textbf{Proof.} We first remark that $span_{A^{\ast}}P(A)=span_{A^{\ast}}P_{1}(A).$ By Corollary 4, $R(A)=\lbrace x\in A: f(x^{\ast}x)=0$ for all $f\in P_{1}(A)\rbrace.$ Let $x\in R(A),$ by the Schwarz inequality, $ \vert f(x)\vert^{2}\leq f(x^{\ast}x)$ for all $f\in P_{1}(A),$ thus $f(x)=0$ for all $f\in P_{1}(A).$ Since  $A'\subset span_{A^{\ast}}P_{1}(A),$ it follows that $g(x)=0$ for all $g\in A',$ hence $x=0.$

\noindent \textbf{}

\noindent \textbf{Remark.} Proposition 8 is a generalization of the following Hennings result: If $A$ is a symmetric unital MQ*-algebra such that $A'\subset span_{A^{\ast}}P(A),$ then $A$ is semisimple.

\noindent \textbf{}

\noindent \textbf{Counter example} Let $E$ be the linear space of absolutely continuous functions on [0,1], with the norm $\Vert f\Vert=\int_{0}^{1}\vert f(t)\vert dt,$ the convolution-multiplication $(f\times g)(x)=\int_{0}^{x}f(x-t)g(t)dt,$ and the natural involution $f^{\ast}(x)=\overline{f(x)}. (E,\Vert .\Vert)$ is a symmetric radical commutative Banach *-algebra [6; p.316]. Let $A=E\bigoplus Ce$ be the algebra obtained from $E$ by adjunction of an identity $e.$ A is a commutative Banach *-algebra with identity, $A$ is symmetric [6; Theorem 4.7.9]. It is easily shown that $R(A)=E,$ thus $A$ is not semisimple, hence $A'$ is not included in $span_{A^{\ast}}P(A)$ by Proposition 8.

\noindent \textbf{}

\noindent \textbf{V. On representations of BP*-algebras}

\noindent \textbf{}

Let $A$ be a BP*-algebra and $f\in P(A).$ Let $L_{f}=\lbrace x\in A: f(x^{\ast}x)=0\rbrace, L_{f}$ is a linear space and the quotient space $A/L_{f}$ is a pre-Hilbert space with the inner product $\langle x+L_{f},y+L_{f}\rangle = f(y^{\ast}x), x,y\in A.$ Since $f$ is hermitian and admissible, we can associate to $f$ a cyclic representation on a Hilbert space $H_{f}$ which is the completion of $A/L_{f}.$ For each  $x\in A, T_{f}x$ is defined by $T_{f}x(y+L_{f})=xy+L_{f}$ for all $y\in A.$

\noindent \textbf{}

\noindent \textbf{Theorem 9.} Let $A$ be a BP*-algebra such that $M_{s}(A)$ is equicontinuous. Then for each $f\in P(A),$ the representation $T_{f}$ is continuous.

\noindent \textbf{}

\noindent \textbf{Proof.} Let $x,y \in A,$ we have $f(y^{\ast}x^{\ast}xy)\leq f(y^{\ast}y)\rho_{s}(x^{\ast}x),$ thus $\langle xy+L_{f},xy+L_{f}\rangle\leq \langle y+L_{f},y+L_{f}\rangle\rho_{s}(x^{\ast}x),$ i.e. $ \Vert T_{f}x(y+L_{f})\Vert\leq\Vert y+L_{f}\Vert\rho_{s}(x^{\ast}x)^{\frac{1}{2}},$
hence $ \Vert T_{f}x\Vert\leq\rho_{s}(x^{\ast}x)^{\frac{1}{2}}=\rho_{s}(x).$ Since $M_{s}(A)$ is equicontinuous, $\rho_{s}$ is continuous, therefore $T_{f}$ is continuous.

\noindent \textbf{}

\noindent \textbf{Remark.} Theorem 9 is an improvement of [5; Theorem 4.2].

\noindent \textbf{}

\noindent \textbf{Acknowledgments.} I would like to express my thanks to Professor T. Husain for help he has given me during the preparation of this work.

\noindent \textbf{}

\begin{center}
\noindent \textbf{REFERENCES}
\end{center}
 
\noindent \textbf{} 
 
\noindent \textbf{} [1] M. Akkar, M. El Azhari, M.Oudadess, On an expression of the spectrum in BP*-algebras, Period. Math. Hungar. 19 (1988), 65-67.

\noindent \textbf{} [2] G. R. Allan, On a class of locally convex algebras, Proc. London Math. Soc.      17 (1967), 91-114.

\noindent \textbf{} [3] M. A. Hennings, Topologies on BP*-algebras, J. Math. Anal. Appl. 140 (1989), 289-300.

\noindent \textbf{} [4] T. Husain, S. A. Warsi, Positive functionals on BP*-algebras, Period. Math. Hungar. 8 (1977), 15-28.

\noindent \textbf{} [5] T. Husain, S. A. Warsi, Representations of BP*-algebras, Math. Japon. 21 (1976), 237-247.

\noindent \textbf{} [6] C. E. Richart, General Theory of Banach Algebras, Van Nostrand, Princeton N. J., 1960. 

\noindent \textbf{} [7] A. P. Robertson, W, Robertson, Topological Vector Spaces. Cambridge Tracts in Math. and Math. Phys. 53, Cambridge Univ. Press, London, 1973.

\noindent \textbf{} 

\noindent \textbf{} Ecole Normale Sup\'{e}rieure

\noindent \textbf{} Avenue Oued Akreuch

\noindent \textbf{} Takaddoum, BP 5118, Rabat

\noindent \textbf{} Morocco
 
\noindent \textbf{} 

\noindent \textbf{} E-mail:  mohammed.elazhari@yahoo.fr

\end{document}